\documentclass[12pt]{article}
\usepackage{latexsym}
\usepackage{amsmath}
\usepackage{amssymb}
\usepackage{amsfonts}

\numberwithin{equation}{section}

\newtheorem{theorem}{Theorem}[section]
\newtheorem{lemma}{Lemma}[section]

\newtheorem{remark}{Remark}[section]
\newcommand{\eproof}{{\mbox{\ }~\hfill
\mbox{\large $\Box$} \par \vskip 10pt}}

\newcommand{\R}{{\mathbb R}}

\title{Quantitative uniqueness for the power of Laplacian with singular coefficients}
\author{Ching-Lung Lin\thanks{Department of Mathematics, National Chung Cheng University,
Chia-Yi 62117, Taiwan. Email:cllin@math.ccu.edu.tw}\qquad Sei
Nagayasu
\thanks{Department of Mathematics, Taida Institute for Mathematical
Sciences, National Taiwan University, Taipei 106, Taiwan.
Email:nagayasu@math.ntu.edu.tw}
 \qquad Jenn-Nan
Wang\thanks{Department of Mathematics, Taida Institute of
Mathematical Sciences, NCTS (Taipei), National Taiwan University,
Taipei 106, Taiwan. Email:jnwang@math.ntu.edu.tw}}

\date{}

\begin{document}
\maketitle

\begin{abstract}
In this paper we study the local behavior of a solution to the $l$th
power of Laplacian with singular coefficients in lower order terms.
We obtain a bound on the vanishing order of the nontrivial solution.
Our proofs use Carleman estimates with carefully chosen weights. We
will derive appropriate three-sphere inequalities and apply them to
obtain doubling inequalities and the maximal vanishing order.
\end{abstract}

\section{Introduction}\label{sec1}
\setcounter{equation}{0}

Assume that $\Omega$ is a connected open set containing $0$ in
$\R^n$ for $n\geq 2$. In this paper we are interested in the local
behavior of $u$ satisfying the following differential inequality:
\begin{equation}\label{1.1}
\begin{array}{rl}
|\triangle^l u| \leq K_0\sum_{|\alpha|\leq
l-1}|x|^{-2l+|\alpha|}|D^\alpha u|+K_0\sum_{|\alpha|=
l}^{[3l/2]}|x|^{-2l+|\alpha|+\epsilon}|D^\alpha u|,
\end{array}
\end{equation}
where $0<\epsilon<1/2$ and $[h]=k\in{\mathbb Z}$ when $k\le h<k+1$.
For \eqref{1.1}, a strong unique continuation was proved by the
first author \cite{Lin}. A similar result for the power of Laplacian
with lower derivatives up to $l$-th order can be found in
\cite{cogr}. On the other hand, a unique continuation property for
the $l$-th power of Laplacian with the same order of lower
derivatives as in \eqref{1.1} was given in \cite{pro}. Those results
mentioned above concern only the qualitative behavior of the
solution. In other words, they showed that if $u$ vanishes at $0$ in
infinite order or $u$ vanishes in an open subset of $\Omega$, then
$u$ must vanishes identically in $\Omega$. The aim of this paper is
to study the strong unique continuation from a quantitative
viewpoint. Namely, we are interested in the maximal vanishing order
at $0$ of any nontrivial solution to \eqref{1.1}. It is worth
mentioning that quantitative estimates of the strong unique
continuation are useful in studying the nodal sets of eigenfunctions
\cite{dofe}, or solutions of second order elliptic equations
\cite{hasi}, \cite{lin2}, or the inverse problem \cite{abrv}.

Perhaps, for the quantitative uniqueness problem, the most popular
technique, introduced by Garofala and Lin \cite{gl1}, \cite{gl2}, is
to use the frequency function related to the solution. This method
works quite efficiently for second order strongly elliptic
operators. However, this method can not be applied to \eqref{1.1}.
Another method to derive quantitative estimates of the strong unique
continuation is based on Carleman estimates, which was first
initiated by Donnelly and Fefferman \cite{dofe} where they studied
the maximal vanishing order of the eigenfunction with respect to the
corresponding eigenvalue on a compact smooth Riemannian manifold.
Their method does not work for \eqref{1.1} either.

Recently, the first and third authors and Nakamura \cite{Lin0}
introduced a method based on appropriate Carleman estimates to prove
a quantitative uniqueness for second order elliptic operators with
sharp singular coefficients in lower order terms. A key strategy of
our method is to derive three-sphere inequalities and then apply
them to obtain doubling inequalities and the maximal vanishing
order. Both steps require delicate choices of cut-off functions.
Nevertheless, this method is quite versatile and can be adopted to
treat many equations or even systems. The present work is an
interesting application of the ideas of \cite{Lin0} to the $l$th
power of Laplacian with singular coefficients. The power $l=2$ is
the most interesting and useful case. It corresponds to the
biharmonic operator with third order derivatives. Our work provides
a quantitative estimate of the strong unique continuation for this
equation. To our best knowledge, this quantitative estimate has not
been derived before.

We now state main results of the paper. Assume that
$B_{R'_0}\subset\Omega$ for some $R'_0>0$.
\begin{theorem}\label{thm1.1}
There exists a positive number $\tilde{R}_0<e^{-1/2}$ such that if
 $\ 0<r_1<r_2<r_3\leq R'_0$ and $r_1/r_3<r_2/r_3<\tilde{R}_0$, then
\begin{equation}\label{1.2}
\int_{|x|<r_2}|u|^2dx\leq
C\left(\int_{|x|<r_1}|u|^2dx\right)^{\tau}\left(\int_{|x|<{r_3}}|u|^2dx\right)^{1-\tau}
\end{equation}
for $u\in H^{2l}({B}_{R'_0})$ satisfying \eqref{1.1} in
${B}_{R'_0}$, where $C$ and $0<\tau<1$ depend on $r_1/r_3$,
$r_2/r_3$, $n$, $l$, and $K_0$.
\end{theorem}
\begin{remark}\label{rem1.0}
From the proof, the constants $C$ and $\tau$ can be explicitly
written as $C=\max\{C_0(r_2/r_1)^n,\exp(B\beta_0)\}$ and
$\tau=B/(A+B)$, where $C_0>1$ and $\beta_0$ are constants depending
on $n,l, K_0$ and
\begin{eqnarray*}
&&A=A(r_1/r_3,r_2/r_3)=(\log (r_1/r_3)-1)^2-(\log (r_2/r_3))^2,\\
&&B=B(r_2/r_3)=-1-2\log (r_2/r_3).
\end{eqnarray*}
The explicit forms of these constants are important in the proof of
Theorem~\ref{thm1.2}.
\end{remark}

\begin{theorem}\label{thm1.2}
Let $u\in H^{2l}_{loc}({\Omega})$ be a nonzero solution to
\eqref{1.1}. Then we can find a constant $R_2$ (depending on $n,
l,\epsilon, K_0$) and a constant $m_1$ (depending on $n, l,
\epsilon,K_0, \|u\|_{L^2(|x|<{R_2^2})}/\|u\|_{L^2(|x|<{R_2^4})}$)
satisfying
\begin{equation}\label{1.3}
\limsup_{R\to 0}\frac{1}{R^{m_1}} \int_{|x|<R}|u|^2 dx>0.
\end{equation}
\end{theorem}

\begin{theorem}\label{thm1.3}
Let $u\in H^{2l}_{loc}({\Omega})$ be a nonzero solution to
\eqref{1.1}. Then there exist positive constants $R_3$ (depending on
$n, l,\epsilon, K_0$) and $C_3$ (depending on $n, l, \epsilon, K_0,
m_1$) such that if $0<r\leq R_3$, then
\begin{equation}\label{1.4}
\int_{|x|\le{2r}}|u|^2dx\leq C_3\int_{|x|\le{r}}|u|^2dx,
\end{equation}
where $m_1$ is the constant obtained in Theorem \ref{thm1.2}.
\end{theorem}

The rest of the paper is devoted to the proofs of
Theorem~\ref{thm1.1}-\ref{thm1.3}.

\section{Three-sphere inequalities}\label{sec2}
\setcounter{equation}{0}

We will prove Theorem~\ref{thm1.1} in this section. To begin, we
recall a Carleman estimate with weight
$\varphi_{\beta}=\varphi_{\beta}(x) =\exp (\frac{\beta}{2}(\log
|x|)^2)$ given in \cite{Lin}.
\begin{lemma}{\rm \cite[Corollary~3.3]{Lin}}\label{lem2.1}
There exist a sufficiently large number $\beta_0>0$ and a
sufficiently small number $r_0>0$, depending on $n$ and $l$, such
that for all $u\in U_{r_0}$ with $0<r_0<e^{-1}$, $\beta\geq
\beta_0$, we have that
\begin{equation}\label{2.1}
\begin{array}{l}
\sum_{|\alpha|\leq 2l}\beta^{3l-2|\alpha|}  \int \varphi^2_\beta
{|x|^{2|\alpha|-n}(\log|x|)^{2l-2|\alpha|}|D^\alpha u|^2
dx}\\
\leq \tilde C_0\int \varphi^2_\beta |x|^{4l-n}|\Delta^l u|^2 dx,
\end{array}
\end{equation}
where $U_{r_0}=\{u\in C_0^{\infty}(\mathbb{R}^n\setminus\{0\}):
\mbox{\rm supp}(u)\subset B_{r_0}\}$ and $\tilde C_0$ is a positive
constant depending on $n$ and $l$. Here $e=\exp(1)$.
\end{lemma}

\begin{remark}\label{rem2.1}
The estimate \eqref{2.1} in Lemma \ref{lem2.1} remains valid if we
assume $u\in H_{loc}^{2l}(\mathbb{R}^n\setminus\{0\})$ with compact
support. This can be easily obtained by cutting off $u$ for small
$|x|$ and regularizing.
\end{remark}

We first consider the case where $0<r_1<r_2<R<1/e$ and
$B_R\subset\Omega$. The constant $R$ will be chosen later. To use
the estimate \eqref{2.1}, we need to cut-off $u$. So let $\xi(x)\in
C^{\infty}_0 ({\mathbb R}^n)$ satisfy $0\le\xi(x)\leq 1$ and
\begin{equation*}
\xi (x)=
\begin{cases}
\begin{array}{l}
0,\quad |x|\leq r_1/e,\\
1,\quad r_1/2<|x|<er_2,\\
0,\quad |x|\geq 3r_2.
\end{array}
\end{cases}
\end{equation*}
It is easy to check that for all multiindex $\alpha$
\begin{equation}\label{gradxi}
\begin{cases}
|D^{\alpha}\xi|=O(r_1^{-|\alpha|})\ \text{for all}\ r_1/e\le |x|\le r_1/2\\
|D^{\alpha}\xi|=O(r_2^{-|\alpha|})\ \text{for all}\ er_2\le |x|\le
3r_2.
\end{cases}
\end{equation}
On the other hand, repeating the proof of Corollary~17.1.4 in
\cite{Hor3}, we can show that
\begin{equation}\label{inter}
\int_{a_1r<|x|<a_2r}||x|^{|\alpha|}D^{\alpha} u|^2dx\le
C'\int_{a_3r<|x|<a_4r}|u|^2dx,\quad|\alpha|\le 2l,
\end{equation}
for all $0<a_3<a_1<a_2<a_4$ such that $B_{a_4r}\subset\Omega$, where
the constant $C'$ is independent of $r$ and $u$.

Noting that the commutator $[\Delta^l,\xi]$ is a $2l-1$ order
differential operator. Applying \eqref{2.1} to $\xi u$ and using
\eqref{1.1}, \eqref{gradxi}, \eqref{inter} implies
\begin{eqnarray}\label{2.2}
&& \sum_{|\alpha|\leq 2l}\beta^{3l-2|\alpha|} \int_{r_1/2<|x|<er_2}
\varphi^2_\beta |x|^{2|\alpha|-n}(\log|x|)^{2l-2|\alpha|}|D^\alpha
u|^2
dx\notag\\
&\leq& \sum_{|\alpha|\leq 2l}\beta^{3l-2|\alpha|}  \int
\varphi^2_\beta
{|x|^{2|\alpha|-n}(\log|x|)^{2l-2|\alpha|}|D^\alpha (\xi u)|^2
dx}\notag\\
&\leq& \tilde C_0\int \varphi^2_\beta |x|^{4l-n}|\Delta^l (\xi u)|^2 dx\notag\\
&\leq& 2\tilde C_0\int\varphi^2_\beta
|x|^{4l-n}\xi^2{(}K_0\sum_{|\alpha|\leq
l-1}|x|^{-2l+|\alpha|}|D^\alpha u|+K_0\sum_{|\alpha|=
l}^{[3l/2]}|x|^{-2l+|\alpha|+\epsilon}|D^\alpha u|{)}^2dx\notag\\
&&+2\tilde C_0\int\varphi^2_\beta
|x|^{4l-n}\big{|}[\Delta^l,\xi]u\big{|}^2
dx\notag\\
&\leq& \tilde C_1\Big{\{}\int_{r_1/2<|x|<er_2} \varphi^2_\beta
(\sum_{|\alpha|\leq l-1}|x|^{2|\alpha|-n}|D^\alpha
u|^2+\sum_{|\alpha|=
l}^{[3l/2]}|x|^{2|\alpha|-n+2\epsilon}|D^\alpha u|^2) dx\notag\\
&&\qquad + \int_{r_1/e<|x|<r_1/2} \varphi^2_\beta
\sum_{|\alpha|\leq
2l-1}|x|^{2|\alpha|-n}|D^\alpha u|^2 dx \notag\\
&&\qquad+\int_{er_2<|x|<3r_2} \varphi^2_\beta \sum_{|\alpha|\leq
2l-1}|x|^{2|\alpha|-n}|D^\alpha u|^2 dx\Big{\}}\notag\\
&\leq& \tilde{C}_2\Big{\{}\int_{r_1/2<|x|<er_2} \varphi^2_\beta
(\sum_{|\alpha|\leq l-1}|x|^{2|\alpha|-n}|D^\alpha
u|^2+\sum_{|\alpha|=
l}^{[3l/2]}|x|^{2|\alpha|-n+2\epsilon}|D^\alpha u|^2) dx\notag\\
&&\qquad + r_1^{-n}\varphi^2_\beta(r_1/e) \int_{r_1/e<|x|<r_1/2}
\sum_{|\alpha|\leq
2l-1}||x|^{|\alpha|}D^{\alpha} u|^2dx\notag\\
&&\qquad + r_2^{-n}\varphi^2_\beta(er_2) \int_{er_2<|x|<3r_2}
\sum_{|\alpha|\leq
2l-1}||x|^{|\alpha|}D^{\alpha} u|^2 dx\Big{\}}\notag\\
&\leq& \tilde{C}_3\Big{\{}\int_{r_1/2<|x|<er_2} \varphi^2_\beta
(\sum_{|\alpha|\leq l-1}|x|^{2|\alpha|-n}|D^\alpha
u|^2+\sum_{|\alpha|=
l}^{[3l/2]}|x|^{2|\alpha|-n+2\epsilon}|D^\alpha u|^2) dx\notag\\
&&\qquad + r_1^{-n}\varphi^2_\beta(r_1/e) \int_{r_1/4<|x|<r_1}
|u|^2dx+r_2^{-n}\varphi^2_\beta(er_2) \int_{2r_2<|x|<4r_2}
|u|^2dx\Big{\}},\notag\\
\end{eqnarray}
where $\tilde C_1$, $\tilde C_2$, and $\tilde C_3$ are independent
of $r_1$, $r_2$, and $u$.

We now choose $r_0<e^{-\epsilon^{-1}([3l/2]-l)-1}$ small enough such
that
$$
\begin{cases}
(\log(er_0))^{-2}\le\frac{1}{2\tilde C_3}\\
(er_0)^{2\epsilon}(\log(er_0))^{2([3l/2]-l)}\le\frac{1}{2\tilde
C_3}.
\end{cases}$$
Letting $R\leq r_0$ and $\beta\geq\beta_0\geq \max\{2\tilde
C_3,1\}$, we can absorb the integral over $r_1/2<|x|<er_2$ on the
right side of \eqref{2.2} into its left side to obtain
\begin{eqnarray}\label{2.3}
&& \int_{r_1/2<|x|<er_2} \varphi^2_\beta (\sum_{|\alpha|\leq
l-1}|x|^{2|\alpha|-n}|D^\alpha u|^2+\sum_{|\alpha|=
l}^{[3l/2]}|x|^{2|\alpha|-n+2\epsilon}|D^\alpha u|^2)
dx\notag\\
&\leq&\tilde C_4\Big{\{} r_1^{-n}\varphi^2_\beta(r_1/e)
\int_{r_1/4<|x|<r_1} |u|^2 dx + r_2^{-n}\varphi^2_\beta(er_2)
\int_{2r_2<|x|<4r_2} |u|^2 dx\Big{\}},\notag\\
\end{eqnarray}
where $\tilde C_4=1/\tilde C_3$. Using \eqref{2.3} we have that
\begin{equation*}
r_2^{-n}\varphi^2_\beta(r_2) \int_{r_1/2<|x|<r_2}|u|^2 dx
\end{equation*}
\begin{eqnarray}\label{2.4}
&\leq&\int_{r_1/2<|x|<er_2} \varphi^2_\beta {|x|^{-n}|u|^2 dx}\notag\\
&\leq&\tilde C_4\Big{\{} r_1^{-n}\varphi^2_\beta(r_1/e)
\int_{r_1/4<|x|<r_1} |u|^2 dx + r_2^{-n}\varphi^2_\beta(er_2)
\int_{2r_2<|x|<4r_2} |u|^2 dx\Big{\}}.\notag\\
\end{eqnarray}
Dividing $r_2^{-n}\varphi^2_\beta(r_2)$ on the both sides of
\eqref{2.4} implies
\begin{eqnarray}\label{2.5}
&&\int_{r_1/2<|x|<r_2}|u|^2 dx\notag\\
&\leq&\tilde C_4\Big{\{}
(r_2/r_1)^{n}[\varphi^2_\beta(r_1/e)/\varphi^2_\beta(r_2)]
\int_{r_1/4<|x|<r_1} |u|^2 dx\notag\\
&&\quad + [\varphi^2_\beta(er_2)/\varphi^2_\beta(r_2)]
\int_{2r_2<|x|<4r_2} |u|^2 dx\Big{\}}\notag\\
&\leq&\tilde C_5\Big{\{}
(r_2/r_1)^{n}[\varphi^2_\beta(r_1/e)/\varphi^2_\beta(r_2)]
\int_{|x|<{r_1}} |u|^2 dx\notag\\
&&\quad +(r_2/r_1)^{n}[\varphi^2_\beta(er_2)/\varphi^2_\beta(r_2)]
\int_{|x|<{4r_2}} |u|^2 dx\Big{\}},
\end{eqnarray}
where $\tilde C_5=\max\{\tilde C_4,1\}$. With such choice of $\tilde
C_5$, we can see that
$$
\tilde
C_5(r_2/r_1)^{n}[\varphi^2_\beta(r_1/e)/\varphi^2_\beta(r_2)]>1
$$
for all $0<r_1<r_2$. Adding $\int_{|x|<{r_1/2}} |u|^2 dx$ to both
sides of \eqref{2.5} and choosing $r_2\leq R=\min\{r_0,1/4\}$, we
get that
\begin{eqnarray}\label{2.6}
&&\int_{|x|<{r_2}}|u|^2 dx\notag\\
&\leq& 2\tilde
C_5(r_2/r_1)^{n}[\varphi^2_\beta(r_1/e)/\varphi^2_\beta(r_2)]
\int_{|x|<{r_1}} |u|^2 dx\notag\\
&&+2\tilde C_5(r_2/r_1)^{n}
[\varphi^2_\beta(er_2)/\varphi^2_\beta(r_2)] \int_{|x|<1} |u|^2 dx.
\end{eqnarray}
By denoting
\begin{eqnarray*}
&&A=\beta^{-1}\,\log[\varphi^2_\beta(r_1/e)/\varphi^2_\beta(r_2)]=(\log r_1-1)^2-(\log r_2)^2>0,\\
&&B=-\beta^{-1}\,\log[\varphi^2_\beta(er_2)/\varphi^2_\beta(r_2)]=-1-2\log
r_2>0,
\end{eqnarray*}
\eqref{2.6} becomes
\begin{eqnarray}\label{22.66}
&&\int_{|x|<{r_2}}|u|^2 dx\notag\\
&\leq& 2\tilde C_5(r_2/r_1)^{n}\Big{\{}\exp(A\beta)\int_{|x|<{r_1}}
|u|^2 dx+\exp(-B\beta) \int_{|x|<1} |u|^2 dx\Big{\}}.\notag\\
\end{eqnarray}

To further simplify the terms on the right hand side of
\eqref{22.66}, we consider two cases. If
$$
\int_{|x|<{r_1}} |u|^2 dx\ne 0
$$
and
$$\exp{(A\beta_0)}\int_{|x|<{r_1}} |u|^2 dx<\exp{(-B\beta_0)}\int_{|x|<{1}} |u|^2 dx,$$
then we can pick a $\beta>\beta_0$ such that
$$
\exp{(A\beta)}\int_{|x|<{r_1}} |u|^2
dx=\exp{(-B\beta)}\int_{|x|<{1}} |u|^2 dx.
$$
Using such $\beta$, we obtain from \eqref{22.66} that
\begin{eqnarray}\label{2.8}
&&\int_{|x|<{r_2}}|u|^2 dx\notag\\
&\leq& 4\tilde C_5(r_2/r_1)^{n}\exp{(A\beta)} \int_{|x|<{r_1}} |u|^2dx\notag\\
&=& 4\tilde
C_5(r_2/r_1)^{n}\left(\int_{|x|<{r_1}}|u|^2dx\right)^{\frac{B}{A+B}}\left(\int_{|x|<{1}}|u|^2dx\right)^{\frac{A}{A+B}}.
\end{eqnarray}
If
$$
\int_{|x|<{r_1}} |u|^2 dx=0,
$$
then it follows from \eqref{22.66} that
$$
\int_{|x|<{r_2}}|u|^2 dx=0
$$ since we can take $\beta$ arbitrarily large. The
three-sphere inequality obviously holds.

On the other hand, if
$$ \exp{(-B\beta_0)}\int_{|x|<{1}} |u|^2dx\leq\exp{(A\beta_0)}\int_{|x|<{r_1}} |u|^2 dx,$$
then we have
\begin{eqnarray}\label{2.9}
&&\int_{|x|<{r_2}}|u|^2 dx\notag\\
&\leq&
\left(\int_{|x|<1}|u|^2dx\right)^{\frac{B}{A+B}}\left(\int_{|x|<1}|u|^2dx\right)^{\frac{A}{A+B}}\notag\\
&\leq&
\exp{(B\beta_0)}\left(\int_{|x|<{r_1}}|u|^2dx\right)^{\frac{B}{A+B}}\left(\int_{|x|<1}|u|^2dx\right)^{\frac{A}{A+B}}.
\end{eqnarray}
Putting together \eqref{2.8}, \eqref{2.9}, and setting $\tilde
C_6=\max\{4\tilde C_5(r_2/r_1)^n,\exp{(B\beta_0)}\}$, we arrive at
\begin{equation}\label{2.99}
\int_{|x|<{r_2}}|u|^2 dx \le \tilde
C_6\left(\int_{|x|<{r_1}}|u|^2dx\right)^{\frac{B}{A+B}}\left(\int_{|x|<1}|u|^2dx\right)^{\frac{A}{A+B}}.
\end{equation}

Now for the general case, we take $\tilde R_0=R$ and consider
$0<r_1<r_2<r_3$ with $r_1/r_3<r_2/r_3\le \tilde R_0$. By scaling,
i.e. defining $\widehat{u}(y):=u(r_3y)$, we derive from \eqref{2.99}
that
\begin{equation}\label{2.10}
\int_{|y|<{r_2/r_3}}|\widehat{u}|^2 dy \leq
C(\int_{|y|<{r_1/r_3}}|\widehat{u}|^2dy)^{\tau}(\int_{|y|<1}|\widehat{u}|^2dy)^{1-\tau},
\end{equation}
where $\tau=B/(A+B)$ with
\begin{eqnarray*}
&&A=A(r_1/r_3,r_2/r_3)=(\log (r_1/r_3)-1)^2-(\log (r_2/r_3))^2,\\
&&B=B(r_2/r_3)=-1-2\log (r_2/r_3),
\end{eqnarray*}
and $C=\max\{4\tilde C_5(r_2/r_1)^n,\exp(B\beta_0)\}$. Note that
that $\tilde C_5$ can be chosen independent of the scaling factor
$r_3$ provided $r_3<1$. Replacing the variable $y=x/r_3$ in
\eqref{2.10} gives
$$
\int_{|x|<{r_2}}|u|^2 dx \leq
C(\int_{|x|<{r_1}}|u|^2dx)^{\tau}(\int_{|x|<{r_3}}|u|^2dx)^{1-\tau}.
$$
This ends the proof. \eproof

\section{Doubling inequalities and maximal vanishing order}\label{sec3}
\setcounter{equation}{0}

In this section, we prove Theorem~\ref{thm1.2} and
Theorem~\ref{thm1.3}. We begin with another Carleman estimate
derived in \cite[Lemma~2.1]{Lin}: for any
 $u\in {C^\infty_0
({\mathbb R}^n \backslash {\{0\}})}$ and for any $m\in
{\{k+1/2,k\in{\mathbb N}\}}$, we have the following estimate
\begin{equation}\label{3.1}
\sum_{|\alpha|\leq 2l} \int m^{2l-2|\alpha|}
|x|^{-2m+2|\alpha|-n}|D^\alpha u|^2 dx\leq  C\int
{|x|^{-2m+4l-n}|\triangle^l u|^2 dx},
\end{equation}
where $C$ depends only on the dimension $n$ and the power $l$.

\begin{remark}\label{rem3.1}
Using the cut-off function and regularization, estimate
\eqref{3.1} remains valid for any fixed $m$ if $u\in
H^{2l}_{loc}({\mathbb R}^n \backslash {\{0\}})$ with compact
support.
\end{remark}

In view of Remark~\ref{rem3.1}, we can apply \eqref{3.1} to the
function $\chi u$ with $\chi(x) \in C^{\infty}_0 ({\mathbb
R}^n\backslash {\{0\}})$. Thus, we define $\chi(x) \in C^{\infty}_0
({\mathbb R}^n\backslash {\{0\}})$ as
$$
\chi(x)=\begin{cases}
0\quad\text{if}\quad
|x|\leq \delta/3,\\
1\quad\text{in}\quad \delta/2\leq|x|\leq(R_0+1)R_0R/4=r_4R,\\
0\quad\text{if}\quad 2r_4R\leq |x| ,
\end{cases}
$$
where $\delta\le R_0^2R/4$, $R_0>0$ is a small number which will be
chosen later and $R<1$ is sufficiently small. Here the number $R$ is
not yet fixed and is given by $R=(\gamma m)^{-l/2\epsilon}$, where
$\gamma>0$ is a large constant which will be determined later. Using
the estimate \eqref{3.1} and the equation \eqref{1.1}, we can derive
that
\begin{eqnarray*}
&& \sum_{|\alpha|\leq 2l} \int_{\delta/2\leq|x|\leq {r_4R}}
m^{2l-2|\alpha|}|x|^{-2m+2|\alpha|-n}|D^\alpha u|^2 dx\\
&\leq&  \sum_{|\alpha|\leq 2l} \int m^{2l-2|\alpha|}
|x|^{-2m+2|\alpha|-n}|D^\alpha (\chi u)|^2 dx\\
&\leq&  C\int |x|^{-2m+4l-n}|\Delta^l (\chi u)|^2 dx\\
&=& C\int_{\delta/2\leq|x|\leq {r_4R}} |x|^{-2m+4l-n}|\Delta^l u|^2
dx+
C\int_{|x|> {r_4R}} |x|^{-2m+4l-n}|\Delta^l (\chi u)|^2 dx\notag\\
&&+C\int_{\delta/3\leq|x|\leq \delta/2} |x|^{-2m+4l-n}|\Delta^l (\chi u)|^2 dx\notag\\
&\leq&  C'K_0^2\int_{\delta/2\leq|x|\leq {r_4R}}
\big{(}\sum_{|\alpha|\leq l-1}|x|^{2|\alpha|-n-2m}|D^\alpha
u|^2+\sum_{|\alpha|=
l}^{[3l/2]}|x|^{2|\alpha|-n-2m+2\epsilon}|D^\alpha u|^2\big{)} dx\notag\\
&&+C\int_{|x|> {r_4R}} |x|^{-2m+4l-n}|\Delta^l(\chi u)|^2 dx\notag
+C\int_{\delta/3\leq|x|\leq \delta/2} |x|^{-2m+4l-n}|\Delta^l(\chi u)|^2 dx\notag\\
\end{eqnarray*}
\begin{eqnarray}\label{3.3}
&\leq&C'K_0^2(r_4R)^{2\epsilon}\int_{\delta/2\leq|x|\leq {r_4R}}
\sum_{|\alpha|=
l}^{[3l/2]}|x|^{2|\alpha|-n-2m}|D^\alpha u|^2 dx\notag\\
&&+C'K_0^2\int_{\delta/2\leq|x|\leq {r_4R}} \sum_{|\alpha|\leq
l-1}|x|^{2|\alpha|-n-2m}|D^\alpha u|^2 dx\notag\\
&&+C\int_{|x|> {r_4R}} |x|^{-2m+4l-n}|\Delta^l(\chi u)|^2 dx
+C\int_{\delta/3\leq|x|\leq \delta/2} |x|^{-2m+4l-n}|\Delta^l(\chi
u)|^2 dx,\notag\\
\end{eqnarray}
where the constant $C'$ depends on $n$ and $l$.

By carefully checking terms on both sides of \eqref{3.3}, we now
choose $\gamma\ge(2C'K_0^2)^{1/l}$ and thus
$$R^{2\epsilon}=(\gamma m)^{-l}\le\frac{m^{-l}}{2C'K_0^2}.$$
Hence, choosing $R_0<1$ (suffices to guarantee
$r_4^{2/\epsilon}=R_0^{2\epsilon}(R_0+1)^{2\epsilon}/4^{2\epsilon}<1$)
and $m$ such that $m^2>2C'K_0^2$, we can remove the first two terms
on the right hand side of the last inequality in \eqref{3.3} and
obtain
\begin{eqnarray}\label{3.31}
&& \sum_{|\alpha|\leq 2l} \int_{\delta/2\leq|x|\leq {r_4R}}
m^{2l-2|\alpha|}|x|^{-2m+2|\alpha|-n}|D^\alpha u|^2 dx\notag\\
&\le& 2C\int_{\delta/3<|x|<\delta/2}|x|^{-2m+4l-n}|\Delta^l (\chi
u)|^2
dx\notag\\
&&+2C\int_{r_4R<|x|<2r_4R}|x|^{-2m+4l-n}|\Delta^l (\chi u)|^2 dx.
\end{eqnarray}

In view of the definition of $\chi$, it is easy to see that for
all multiindex $\alpha$
\begin{equation}\label{chii}
\begin{cases}
|D^{\alpha}\chi|=O(\delta^{-|\alpha|})\ \text{for all}\ \delta/3<|x|<\delta/2,\\
|D^{\alpha}\chi|=O((r_4R)^{-|\alpha|})\ \text{for all}\
r_4R<|x|<2r_4R.
\end{cases}
\end{equation}
Note that $R_0^2\le r_4$ provided $R_0\le 1/3$. Therefore, using
\eqref{chii} and \eqref{inter} in \eqref{3.31}, we derive
\begin{eqnarray}\label{3.4}
&&m^2(2\delta)^{-2m-n}\int_{\delta/2<|x|\le
2\delta}|u|^2dx+m^2(R_0^2R)^{-2m-n}\int_{2\delta<
|x|\le R_0^2R}|u|^2dx\notag\\
&\le&  \sum_{|\alpha|\leq 2l} \int_{\delta/2\leq|x|\leq {r_4R}}
m^{2l-2|\alpha|}|x|^{-2m+2|\alpha|-n}|D^\alpha u|^2 dx\notag\\
&\le& C''\sum_{|\alpha|\le
2l}\delta^{-4l+2|\alpha|}\int_{\delta/3<|x|<\delta/2}|x|^{-2m+4l-n}|D^{\alpha}u|^2 dx\notag\\
&&+C''\sum_{|\alpha|\le
2l}(r_4R)^{-4l+2|\alpha|}\int_{r_4R<|x|<2r_4R}|x|^{-2m+4l-n}|D^{\alpha}u|^2
dx\notag\\
&\le&\tilde C'\delta^{-2m-n}\int_{|x|\le\delta}|u|^2
dx+C''(r_4R)^{-2m-n}\int_{|x|\le R_0R}|u|^2 dx,
\end{eqnarray}
where $\tilde C'=C''3^{2m+n}$ and $C''$ is independent of $R_0$,
$R$, and $m$.

We then add $m^2(2\delta)^{-2m-n}\int_{|x|\le\delta/2}|u|^2dx$ to
both sides of \eqref{3.4} and obtain
\begin{eqnarray}\label{3.5}
&&\frac{1}{2}m^2(2\delta)^{-2m-n}\int_{|x|\le 2\delta}|u|^2dx+m^2(R_0^2R)^{-2m-n}\int_{|x|\le R_0^2R}|u|^2dx\notag\\
&=&\frac{1}{2}m^2(2\delta)^{-2m-n}\int_{|x|\le 2\delta}|u|^2dx+m^2(R_0^2R)^{-2m-n}\int_{|x|\le 2\delta}|u|^2dx\notag\\
&&+m^2(R_0^2R)^{-2m-n}\int_{2\delta<|x|\le R_0^2R}|u|^2dx\notag\\
&\le&\frac{1}{2}m^2(2\delta)^{-2m-n}\int_{|x|\le 2\delta}|u|^2dx+\frac{1}{2}m^2(2\delta)^{-2m-n}\int_{|x|\le 2\delta}|u|^2dx\notag\\
&&+m^2(R_0^2R)^{-2m-n}\int_{2\delta<|x|\le R_0^2R}|u|^2dx\notag\\
&\le&\tilde
C''\delta^{-2m-n}\int_{|x|\le\delta}|u|^2dx+C''(r_4R)^{-2m-n}\int_{|x|\le
R_0R}|u|^2 dx\notag\\
&=&\tilde C''\delta^{-2m-n}\int_{|x|\le\delta}|u|^2dx\notag\\
&&+m^2(R_0^2R)^{-2m-n}C''m^{-2}(\frac{R_0^2}{r_4})^{2m+n}\int_{|x|\le
R_0R}|u|^2 dx
\end{eqnarray}
with $\tilde C''=\tilde C'+2^{2m+n}m^2$.

We first observe that
\begin{eqnarray*}
&&
C''m^{-2}(\frac{R_0^2}{r_4})^{2m+n}=C''m^{-2}\left(\frac{4R_0}{R_0+1}\right)^{2m+n}\notag\\
&\le&C''m^{-2}(4R_0)^{2m+n}\ \le \exp(-2m)
\end{eqnarray*}
for all $R_0\le 1/16$ and $m^2\ge C''$. Thus, we obtain that
\begin{eqnarray}\label{3.6}
&&\frac{1}{2}m^2(2\delta)^{-2m-n}\int_{|x|\le 2\delta}|u|^2dx+m^2(R_0^2R)^{-2m-n}\int_{|x|\le R_0^2R}|u|^2dx\notag\\
&\leq&\tilde C''\delta^{-2m-n}\int_{|x|\le\delta}|u|^2dx\notag\\
&&+m^2(R_0^2R)^{-2m-n}\exp(-2m)\int_{|x|\le R_0R}|u|^2 dx.
\end{eqnarray}

It should be noted that \eqref{3.6} is valid for all $m=j+\frac 12$
with $j\in{\mathbb N}$ and $j\ge j_0$, where $j_0$ depends on $n$,
$l$, $\epsilon$, and $K_0$. Setting $R_j=(\gamma(j+\frac
12))^{-l/2\epsilon}$ and using the relation $m=(\gamma
)^{-1}(R)^{-2\epsilon/l}$, we get from \eqref{3.6} that
\begin{eqnarray}\label{3.7}
&&\frac{1}{2}m^2(2\delta)^{-2m-n}\int_{|x|\le 2\delta}|u|^2dx+m^2(R_0^2R_j)^{-2m-n}\int_{|x|\le R_0^2R_j}|u|^2dx\notag\\
&\leq&\tilde
C''\delta^{-2m-n}\int_{|x|\le\delta}|u|^2dx\notag\\
&&+m^2(R_0^2R_j)^{-2m-n}\exp(-2cR_j^{-2\epsilon/l})\int_{|x|\le
R_0R_j}|u|^2 dx
\end{eqnarray}
for all $j\ge j_0$ and $c=\gamma^{-1}$. We now let $j_0$ be large
enough such that
$$
R_{j+1}<R_j<2R_{j+1}\quad\text{for all}\quad j\ge j_0.
$$
Thus, if $R_{j+1}<R\leq R_j$ for $j\ge j_0$, we can conclude that
\begin{eqnarray}\label{3.8}
\begin{cases}
&\int_{|x|\leq R_0^2R} |u|^2 dx
\leq \int_{|x|\leq R_0^2R_j} |u|^2 dx,\\
&\exp(-2cR_j^{-2\epsilon/l})\int_{|x|\leq R_0R_j} |u|^2 dx
\leq \exp(-cR^{-2\epsilon/l})\int_{|x|\leq R} |u|^2 dx,
\end{cases}
\end{eqnarray}
where we have used the inequality $R_0R_j\le R_{j}/16<R_{j+1}$ to
derive the second inequality above. Namely, we have from
\eqref{3.7} and \eqref{3.8} that
\begin{eqnarray}\label{3.9}
&&\frac{1}{2}m^2(2\delta)^{-2m-n}\int_{|x|\le 2\delta}|u|^2dx
+m^2(R_0^2R_j)^{-2m-n}\int_{|x|\le R_0^2R}|u|^2dx\notag\\
&\leq&\tilde
C''\delta^{-2m-n}\int_{|x|\le\delta}|u|^2dx\notag\\
&&+m^2(R_0^2R_j)^{-2m-n}\exp(-cR^{-2\epsilon/l})\int_{|x|\le
R}|u|^2 dx.
\end{eqnarray}

If there exists $s\in{\mathbb N}$ such that
\begin{equation}\label{3.10}
R_{j+1}<R_0^{2s}\le R_j\quad\text{for some}\quad j\ge j_0,
\end{equation}
then replacing $R$ by $R_0^{2s}$ in \eqref{3.9} leads to
\begin{eqnarray}\label{3.11}
&&\frac{1}{2}m^2(2\delta)^{-2m-n}\int_{|x|\le 2\delta}|u|^2dx
+m^2(R_0^2R_j)^{-2m-n}\int_{|x|\le R_0^{2s+2}}|u|^2dx\notag\\
&\leq&\tilde
C''\delta^{-2m-n}\int_{|x|\le\delta}|u|^2dx\notag\\
&&+m^2(R_0^2R_j)^{-2m-n}\exp(-cR_0^{-4s\epsilon/l})\int_{|x|\le
R_0^{2s}}|u|^2 dx.
\end{eqnarray}
Here $s$ and $R_0$ are yet to be determined. The trick now is to
find suitable $s$ and $R_0$ satisfying \eqref{3.10} and the
inequality
\begin{equation}\label{3.21}
\exp(-cR_0^{-4s\epsilon/l})\int_{|x|\leq R_0^{2s}} |u|^2
dx\leq\frac{1}{2}\int_{|x|\leq R_0^{2s+2}} |u|^2 dx
\end{equation}
holds with such choices of $s$ and $R_0$.

It is time to use the three-sphere inequality \eqref{1.2}. To this
end, we choose $r_1=R_0^{2k+2}$, $r_2=R_0^{2k}$ and $r_3=R_0^{2k-2}$
for $k\ge 1$. Note that $r_1/r_3<r_2/r_3\le R_0^{2}\le \tilde{R}_0$.
Thus \eqref{1.2} implies
\begin{equation}\label{3.12}
\int_{|x|<R_0^{2k}}|u|^2dx/\int_{|x|<R_0^{2k+2}}|u|^2dx\leq
C^{1/\tau}(\int_{|x|<R_0^{2k-2}}|u|^2dx/\int_{|x|<R_0^{2k}}|u|^2dx)^{a},
\end{equation}
where
$$C=\max\{C_0R_0^{-2n},\exp(\beta_0(-1-4\log R_0))\}$$ and
\begin{eqnarray*}
a=\frac{1-\tau}{\tau}=\frac{A}{B}&=&\frac{(\log (r_1/r_3)-1)^2-(\log
(r_2/r_3))^2}{-1-2\log (r_2/r_3)}\\
&=&\frac{(4\log R_0-1)^2-(2\log R_0)^2}{-1-4\log R_0}.
\end{eqnarray*}
It is not hard to see that
\begin{equation}\label{3.13}
\begin{cases}
1<C\le C_0 R_0^{-\beta_1},\\
2< a\le -4\log R_0,
\end{cases}
\end{equation}
where $\beta_1=\max\{2n,4\beta_0\}$ and if $R_0$ is sufficiently
small, e.g., $R_0\le e^{-4}$. Combining \eqref{3.13} and using
\eqref{3.12} recursively, we have that
\begin{eqnarray}\label{3.14}
&&\int_{|x|\leq R_0^{2s}} |u|^2 dx/\int_{|x|\leq R_0^{2s+2}} |u|^2
dx\notag\\
&\leq&C^{1/\tau}(\int_{|x|<R_0^{2s-2}}|u|^2dx/\int_{|x|<R_0^{2s}}|u|^2dx)^{a}\notag\\
&\leq&
C^{\frac{a^{s-1}-1}{\tau(a-1)}}(\int_{|x|<R_0^{2}}|u|^2dx/\int_{|x|<R_0^{4}}|u|^2dx)^{a^{s-1}}
\end{eqnarray}
for all $s\ge 1$. Now from the definition of $a$, we have
$\tau=1/(a+1)$ and thus
$$
\frac{a^{s-1}-1}{\tau(a-1)}=\frac{a+1}{a-1}(a^{s-1}-1)\le
3a^{s-1}.
$$
Then it follows from \eqref{3.14} that
\begin{eqnarray}\label{3.15}
&&\int_{|x|\leq R_0^{2s}} |u|^2 dx/\int_{|x|\leq R_0^{2s+2}} |u|^2
dx\notag\\
&\leq& C^{3(-4\log
R_0)^{s-1}}(\int_{|x|<R_0^{2}}|u|^2dx/\int_{|x|<R_0^{4}}|u|^2dx)^{a^{s-1}}\notag\\
&\leq& (C_0^3(R_0)^{-3\beta_1})^{(-4\log
R_0)^{s-1}}(\int_{|x|<R_0^{2}}|u|^2dx/\int_{|x|<R_0^{4}}|u|^2dx)^{a^{s-1}}.
\end{eqnarray}
Thus, by \eqref{3.15}, we can get that
\begin{eqnarray}\label{3.16}
&&\exp(-cR_0^{-4s\epsilon/l})\int_{|x|\leq R_0^{2s}} |u|^2 dx\notag\\
&\leq&
\exp(-cR_0^{-4s\epsilon/l})(C_0^3(R_0)^{-3\beta_1})^{(-4\log
R_0)^{s-1}}\notag\\
&&(\int_{|x|<R_0^{2}}|u|^2dx/\int_{|x|<R_0^{4}}|u|^2dx)^{a^{s-1}}\int_{|x|\leq
R_0^{2s+2}} |u|^2 dx.\notag\\
\end{eqnarray}

Let $\mu=-\log R_0$, then if $R_0\ (\le
\min\{e^{-4},\sqrt{\tilde{R_0}}\})$ is sufficiently small, i.e.,
$\mu$ is sufficiently large, we can see that
$$
4t\epsilon\mu/l >(t-1)\log(4\mu)+\log(\log
C_0^3+3\beta_1\mu)-\log(c/4),
$$
for all $t\in{\mathbb N}$. In other words, we have that for $R_0$
small
\begin{equation}\label{3.17}
(C_0^3R_0^{-3\beta_1})^{(-4\log
R_0)^{t-1}}<\exp(cR_0^{-4t\epsilon/l}/4)<(1/2)\exp(cR_0^{-4t\epsilon/l}/2),
\end{equation}
for all $t\in{\mathbb N}$. We now fix such $R_0$ so that
\eqref{3.17} holds and $$-\frac{4\varepsilon}{l}\log R_0-2\log
a>0.$$ It is a key step in our proof that we can find a universal
constant $R_0$. After fixing $R_0$, we then define a number $t_0$,
depending on $R_0$ and $u$, as
\begin{eqnarray*}
t_0&=&(\log
2-\log(ac)+\log\log(\int_{|x|<R_0^{2}}|u|^2dx/\int_{|x|<R_0^{4}}|u|^2dx))\\
&&\qquad\times(-\frac{4\varepsilon}{l}\log R_0-\log a)^{-1}.
\end{eqnarray*}
With the choice of $t_0$, we can see that
\begin{equation}\label{3.20}
(\int_{|x|<R_0^{2}}|u|^2dx/\int_{|x|<R_0^{4}}|u|^2dx)^{a^{t-1}}\le\exp(cR_0^{-4t\epsilon/l}/2)
\end{equation}
for all $t\ge t_0$.

Let $s_1$ be the smallest positive integer such that $s_1\geq t_0$.
If
\begin{equation}\label{3.18}
R_0^{2s_1}\le R_{j_0}=(\gamma(j_0+1/2))^{-l/2\epsilon},
\end{equation}
then we can find a $j_1\in{\mathbb N}$ with $j_1\ge j_0$ such that
\eqref{3.10} holds, i.e.,
$$
R_{j_1+1}<R_0^{2s_1}\le R_{j_1}.
$$
On the other hand, if
\begin{equation}\label{3.19}
R_0^{2s_1}> R_{j_0},
\end{equation}
then we pick the smallest positive integer $s_2>s_1$ such that
$R_0^{2s_2}\le R_{j_0}$ and thus we can also find a
$j_1\in{\mathbb N}$ with $j_1\ge j_0$ for which \eqref{3.10}
holds. We now define
\begin{equation*}
s=\begin{cases} s_1\quad\text{if}\quad\eqref{3.18}\quad\text{holds},\\
s_2\quad\text{if}\quad\eqref{3.19}\quad\text{holds}.
\end{cases}
\end{equation*}
It is important to note that with such $s$, \eqref{3.10} is
satisfied for some $j_1$ and \eqref{3.17}, \eqref{3.20} hold.
Therefore, we set $m_1=n+2(j_1+1/2)$ and $m=(m_1-n)/2$. Combining
\eqref{3.16}, \eqref{3.17} and \eqref{3.20} yields that
\begin{eqnarray*}
&&\exp(-cR_0^{-4s\epsilon/l})\int_{|x|\leq R_0^{2s}} |u|^2 dx\notag\\
&\leq&
\exp(-cR_0^{-4s\epsilon/l})(C_0^3(R_0)^{-3\beta_1})^{(-3\log
R_0)^{s-1}}\notag\\
&&(\int_{|x|<R_0^{2}}|u|^2dx/\int_{|x|<R_0^{4}}|u|^2dx)^{a^{(s-1)}}\int_{|x|\leq
R_0^{2s+2}} |u|^2 dx.\notag\\
&\leq&\frac{1}{2}\int_{|x|\leq R_0^{2s+2}} |u|^2 dx
\end{eqnarray*}
which is \eqref{3.21}. Using \eqref{3.21} in \eqref{3.11}, we have
that
\begin{eqnarray}\label{3.22}
&&\frac{1}{2}m^2(2\delta)^{-2m-n}\int_{|x|\le 2\delta}|u|^2dx
+\frac{1}{2}m^2(R_0^2R_{j_1})^{-2m-n}\int_{|x|\le R_0^{2s+2}}|u|^2dx\notag\\
&\leq&\tilde C''\delta^{-2m-n}\int_{|x|\le\delta}|u|^2dx.
\end{eqnarray}
From \eqref{3.22}, we get that
\begin{eqnarray}\label{3.23}
\frac{(m_1-n)^2}{8\tilde C''}(R_0^2R_{j_1})^{-m_1}\int_{|x|\le
R_0^{2s+2}}|u|^2dx\leq\delta^{-m_1}\int_{|x|\le\delta}|u|^2dx
\end{eqnarray}
and
\begin{eqnarray*}
\frac{1}{2}m^2(2\delta)^{-2m-n}\int_{|x|\le
2\delta}|u|^2dx\leq\tilde
C''\delta^{-2m-n}\int_{|x|\le\delta}|u|^2dx
\end{eqnarray*}
which implies
\begin{eqnarray}\label{3.24}
\int_{|x|\le 2\delta}|u|^2dx\leq\frac{8\tilde
C''}{(m_1-n)^2}2^{m_1}\int_{|x|\le\delta}|u|^2dx.
\end{eqnarray}

The estimates \eqref{3.23} and \eqref{3.24} are valid for all
$\delta\le R_0^{2s+2}/4$. Therefore, \eqref{1.3} holds with
$R_2=R_0$.  \eqref{1.4} holds with $R_3=R_0^{2s+2}/8$ and
$C_3=\frac{8\tilde C''}{(m_1-n)^2}2^{m_1}$ and the proof is now
complete.\eproof

\section*{Acknowledgements}

The first and third authors are supported in part by the National
Science Council of Taiwan.

\end{document}